\documentclass[12pt,draft]{article}
\usepackage{amssymb,latexsym,amsmath,euscript,mathtools}
\usepackage[T2C]{fontenc}
\usepackage[russian]{babel}
\usepackage{MnSymbol, mathrsfs}

\begin{document}

\title{Approximation Properties of Locally Convex Spaces and the Problem of Uniqueness of the Trace of a Linear Operator \thanks{\* Originally published in Teor. Funktii Funktsional. Anal. i  Prilozhen. No. 39 (1983), 73-87. Original article submitted January 4, 1981. Translated by Robert H. Silverman}\*}
\author{G. L. Litvinov}
\date{}
\maketitle

\begin{tabbing}
\hspace{14.0cm}\=\kill
Selecta Mathematica Sovietica\>
0272-9903/92/010025-16 \$1.50 + 0.20/0\'\\
Vol. 11, No. 1 (1992), p. 25-40.\> \copyright 1992 Birkh\"{a}user Verlag, Basel\'
\end{tabbing}

\section*{\small 1. Introduction}

Linear positive operators in Hilbert spaces that have a trace were described for the first time in studies by von Neumann devoted to the mathematical 
apparatus of quantum mechanics. Subsequently, Schatten and von Neumann [{\bfseries 13}] described all the linear continuous operators that have a trace in Hilbert 
spaces. An operator that possesses a trace was shown to correspond in one-to-one fashion with a bivalent tensor, and it was proved that the trace of an 
operator could be defined as the result of a total contraction of the corresponding tensor. This construction was generalized by Ruston to the case of Banach 
spaces and, independently of Ruston, Grothendieck [{\bfseries 4}] genera-lized the construction to the case of arbitrary locally convex spaces. It was shown that given an 
inductive topological tensor product $E' \bar{\otimes} E$ of some space $E$ and the space $E'$ dual to it, a linear operator in $E$ could be associated with every element of 
this product. However, the \textquotedblleft operators with trace''   obtained in this way are not, generally speaking, continuous; a narrower class is formed by the Fredholm 
operators in the sense of [{\bfseries 4}], which are always weakly continuous. Finally, an even narrower class is formed by the nuclear operators (which are continuous). 
For Banach spaces, these three classes of operators coincide.

Unfortunately, even in a Banach space the same operator may be associa-ted with different tensors. In this case, the trace cannot be defined unambigu-ously for 
an arbitrary nuclear operator. If the trace of an arbitrary Fredholm operator in a locally convex space $E$ is well-defined, it is said that for this space the 
{\it uniqueness problem} has a positive solution [{\bfseries 2}]. The uniqueness problem is closely related to the well-known approximation problem, and also to the basis 
problem, the statement of which goes back to Banach [{\bfseries 2}], [{\bfseries 4}].

Grothendieck proved [{\bfseries 4}] that the uniqueness problem for a Banach space $E$ has a positive solution if and only if $E$ possesses the approximation prope-rty. For locally convex spaces, Grothendieck proved the following \textquotedblleft conditional assertion'' [{\bfseries 4}]: In every locally convex space, all Fredholm operators have a well-defined trace if and only if the approximation problem has a positive solution (i.e., if and only if all Banach spaces have the approximation prope-rty). 
Enflo's well-known counterexample [{\bfseries 2}] to the approximation problem made this \textquotedblleft conditional assertion''  vacuous and demonstrated the existence of a reflexive Banach space in which there is no well-defined trace for every nuclear operator.

In the present article, it is proved that every nuclear operator in a locally convex space $E$ has a well-defined trace if $E$ possesses the approximation 
property. However, even if a space possesses the approximation property this still does not guarantee a positive solution of the uniqueness problem. Below, 
we present an example of a quasi-complete space with the approximation property in which it is not possible to define the trace for all Fredholm operators. 
We prove that the uniqueness problem has a positive solution if $E$ possesses the \textquotedblleft bounded approximation property.\textquotedblright \   Preliminary information and results are presented in \S2. A number of approximation-type properties of locally convex spaces and relations between these properties are considered in \S3. The principal results of the present study, along with certain corollaries from these results (for example, the existence of a matrix trace), may be found in \S4. The present paper extends some of our previous studies [{\bfseries 7}]-[{\bfseries 9}] and includes some of the results referred to in [{\bfseries 9}]; other results from [{\bfseries 9}] are presented in detail in [{\bfseries 8}].

\section*{\small 2. Fredholm operators and nuclear operators. The trace of an operator}

{\bf 2.1.} Suppose that $E$ is a separated (Hausdorff) locally convex space over the field of real or complex numbers\footnote{Below, unless otherwise stipulated, it will be assumed that all the locally convex spaces are Hausdorff and are defined over the same field (of real or complex numbers).}.  We denote by $E'$ the space that is strongly dual to $E$, and by $\langle x', x\rangle$ the value of the functional $x' \in E'$ on the element $x \in E$. Suppose that $Â$ is an absolutely convex set in $E$ (i.e., $Â$ is convex and         $\lambda B \subset B$ for all numbers $\lambda$ such that $\vert \lambda \vert \leq 1$). If $B$ is a bounded set, we denote by $E_B$ the linear (not necessarily closed) subspace in $E$ generated by the elements of $B$ and equipped with the norm defined by $\| x \| =$ inf $\{\lambda  > 0 \colon x \in \lambda B\}$. If $B$ is a complete set (in particular, if it is compact or weakly compact), the normed space $E_B$ is complete, i.e., it is a Banach space. The canonical embedding $E_B \to E$ is continuous.

Suppose that $U$ is an absolutely convex neighborhood of zero in $E$ and that $p_U (x)$ is the corresponding Minkowski functional, i.e., $p_U(u) \mathrel{\mathop:}=$ inf $\{\lambda > 0 \colon x \in \lambda U)$. We let $\widehat{E}_U$ denote the Banach space that is obtained from $E$ as a result of factorization relative to the linear subspace $\{x \in E \colon p_U (x) = 0\}$ and completed with respect to the norm generated by $p_U$. The canonical mapping is continuous.

Let $U^0$ denote the absolute polar of the neighborhood $U$, i.e., $U^0 = \{x' \in E' \colon \vert \langle x', x \rangle \vert \leqslant 1\}$ for all $x \in U$. The set $U^0$ is absolutely convex, bounded, closed, and complete. It is easily seen that the Banach space $E'_{U^0}$ may be considered as a strong dual to $E_U$.

Banach spaces of type $E_B$ and $\widehat{E}_U$ have been systematically treated by Grothendieck [{\bfseries 4}].
\\

{\bf 2.2.}    A linear operator $A$ that maps the space $E$ into itself is called a {\it Fredholm operator} if it may be represented in the form

\begin{equation}
A \colon x \mapsto \sum_{i=1}^{\infty}\lambda_i \langle x'_i, x \rangle x_i ,
\end{equation}
where $\{\lambda_i\}$ is a summable sequence of numbers (i.e., $\sum_{i=1}^\infty \vert \lambda_i \vert < \infty$) àïd the sequences $\{x_i\}$ and $\{x'_i\}$ are contained in the absolutely convex sets $B \subset E$ and $B' \subset E'$, respectively, such that $E_B$ and $E'_{B'}$ are complete spaces.

If the sequence $\{x'_i\}$ is equicontinuous, i.e., if it is possible to set $B' = U^0$ for some neighborhood $U$ of the zero in $E$, the Fredholm operator $A$ is called a 
{\it nuclear operator}.

Every Fredholm operator is weakly continuous and is a nuclear operator relative to the Mackey topology $\tau(E, E')$, i.e., the strongest locally convex topology 
in $E$ such that the dual space coincides with $E'$. Every nuclear operator is a Fredholm operator and is completely continuous (compact). Clearly, if $E$ is a 
Mackey space, i.e., if the topology in $E$ coincides with $\tau(E, E')$ (in particular, if $E$ is a barrelled space), the class of Fredholm operators coincides with 
the class of nuclear operators.

If $E$ is a quasi-complete barrelled space (for example, a Banach or comp-lete metrizable or reflexive space), then $A$ will be a nuclear operator if and only if 
it can be represented in the form (1), where the sequences $\{x_i\}$ and $\{x'_i\}$ are bounded and the sequence $\{\lambda_i\}$ of numbers is summable.
\\

{\bf 2.3.} The mapping $x \mapsto \lambda_i \langle x', x\rangle $ is an operator of rank 1, so that any Fredholm operator may be decomposed into a sum of one-dimensional opera-tors. Every continuous linear operator of finite rank (i.e., an operator with finite-dimensional image) has a decomposition into a finite sum of one-dimen-sional 
operators

\begin{equation}
A \colon x \mapsto \sum_{i=1}^n\lambda_i \langle x'_i, x \rangle x_i ,
\end{equation}
where all coefficients $\lambda_i$, may be set equal to 1. We denote by sp$A$ the trace of the finite-dimensional operator that $A$ induces in its image. It is easy to 
show that sp$A = \sum_{i=1}^n \lambda_i \langle x_i', x_i \rangle $. It is clear that this quantity is independent of the decomposition (2).

If $A$ is a Fredholm operator that can be represented in the form (1), then the series

\begin{equation}
\sum_{i=1}^{\infty}\lambda_i \langle x'_i, x_i \rangle
\end{equation}
converges absolutely, since the quantity $\langle x'_i, x_i \rangle $ is bounded in $i$. Grothendieck [{\bfseries 4}] called the sum of this series the {\it trace} tr$A$ {\it of the Fredholm operator A} under the condition that this quantity is independent of the decomposition (1) of $A$ into a sum of operators of rank 1. In this case, the trace of the operator is {\it well-defined}. Analogously, the {\it trace of the nuclear operator $A$ is well-defined} if the sum of the series (3) is independent of the decomposition 
of $A$ into a sum of the form (1) with the additional condition that the sequence $\{ x'_i \}$ be equicontinuous.
\\

{\bf 2.4.} Let $F$ and $G$ be locally convex spaces. In the algebraic tensor product $F \otimes G$ we consider the strongest locally convex topology such that the canoni-cal bilinear mapping $F \times G \to F \otimes G$ is separately continuous (under this   mapping, the pair $(f, g)$, where  $f \in  F$ and  $g \in  G$,  is  carried  into  the  element $f \otimes g \in F \otimes G$). The completion of $F \otimes G$ in this topology, called the {\it inductive topology}, is denoted by $F \bar{\otimes} G$. The space $F \bar{\otimes} G$ is called the {\it completed induc-tive topological 
tensor product of the spaces $F$ and $G$} [{\bfseries 4}]. Suppose that $H$ is any complete locally convex space. Associating with every linear continuous mapping $F \otimes G \to H$ its composition with the canonical mapping $F \times G \to F \bar{\otimes} G$ yields a one-to-one correspondence between the continuous linear map-pings $F \bar{\otimes} G \to H$ and the separately continuous bilinear mappings $F \times G \to H$.

The mapping $(x', x) \mapsto \langle x', x \rangle $ is a separately continuous bilinear form on $E' \times E$, so that this form may be continued to a continuous linear form on $E' \bar\otimes E$, called the {\it trace}. The trace of the tensor $u \in E' \bar\otimes E$ is denoted by tr $u$. Under the natural correspondence between tensors and linear operators to be described below, the trace of a tensor corresponds to the trace of a linear operator.
\\

{\bf 2.5.} Let $S(E)$ be the algebra of all weakly continuous linear operators in the space $E$ equipped with the weak operator topology. It is specified by the set of 
semi-norms $A \mapsto \vert \langle x', Ax\rangle \vert$, where $A \in S(E)$ and $x'$ and $x$ run through the spaces $E'$ and $E$, respectively.

We denote by $\Gamma$ the linear mapping $E' \otimes E \to S(E)$, under which the element

$$
u = \sum_{i=1}^n \lambda_i  x'_i \otimes x_i
$$
is carried into the finite-rank operator

$$
x \mapsto \sum_{i=1}^n \lambda_i \langle x'_i, x \rangle x_i .
$$

The mapping $\Gamma$ establishes a one-to-one correspondence between $E' \otimes E$ and the space of continuous operators of finite rank; furthermore, tr $u =$ sp $\Gamma (u)$, i.e., the trace of the tensor $u \in E' \otimes E$ coincides with the trace of the operator induced by $\Gamma (u)$ in its finite-dimensional image.

The mapping $\Gamma$ is continuous in the inductive topology, but it is not always possible to continue it to $E' \bar{\otimes} E$, since $S(E)$ is not, generally speaking, a complete space. Grothendieck [{\bfseries 4}] found a linear subspace in $E' \bar{\otimes} E$ to which it is always possible to continue the mapping $\Gamma$. This subspace, denoted by $E' \overbracket[1pt][1pt]{\otimes} E$, consists of all elements of the form

\begin{equation}
u = \sum_{i=0}^\infty \lambda_i  x'_i \otimes x_i,
\end{equation}
where the sequences $\{\lambda_i\}$, $\{ x'_i \}$, and $\{ x_i \}$ are the same as in the decompo-sition (1). In this case, the series (4) converges in $E' \bar{\otimes} E$, and its sum is called the {\it Fredholm kernel}. The continuation of $\Gamma$ onto $E' \overbracket[1pt][1pt]{\otimes} E$ or onto $E' \bar{\otimes} E$ is denoted again by $\Gamma$.

It is clear that if a tensor $u$ is a Fredholm kernel with decomposition (4), then $\Gamma(u)$ is a Fredholm operator with decomposition (1). Thus, the Fredholm 
operators coincide with the images of the Fredholm kernels under the mapping $\Gamma$. The next result follows directly from the definition of the trace of a tensor 
and that of a Fredholm operator: If $\Gamma (u) = 0$ implies tr$u = 0$, then the trace tr$A$ of the Fredholm operator $A$ is well-defined and tr$A =$ tr$\Gamma^{-1}(A)$. This condition is satisfied if the mapping $\Gamma$ is one-to-one. If this condition is satisfied for a linear subspace in $E' \overbracket[1pt][1pt]{\otimes} E$ , then the trace is well-defined for operators that form a narrower class than the class of Fredholm operators (for example, for nuclear operators). For every locally convex space $E$, 
Grothendieck [{\bfseries 4}] has found a linear subspace in $E' \overbracket[1pt][1pt]{\otimes} E$ on which the mapping $\Gamma$ is one-to-one. This subspace consists of tensors that can be represented in the form (4) with the additional condition

$$
\sum_{i=1}^{\infty} \vert \lambda_i \vert^{2/3} < \infty .
$$
In the literature, other classes of operators in Banach spaces have been described (which are narrower than the class of nuclear operators) in which the 
trace is well-defined.

Under certain additional conditions on the space $E$ (for example, if $E$ is complete and metrizable), the mapping $\Gamma$ may be continued to the completed inductive 
topological tensor product $E' \bar{\otimes} E$. Grothendieck [{\bfseries 4}] called operators of the form $\Gamma(u)$, where $u \in E' \bar{\otimes} E$, {\it operators with trace}. Virtually nothing is known about operators with trace that are not Fredholm operators other than the fact that they exist. If $E$ is a Banach space, $E' \bar{\otimes} E$ contains only Fredholm kernels, so that in this case the class of operators with trace, the class of Fredholm operators, and the class of nuclear operators all coincide.
\\

{\bfseries Remark.} The class of nuclear operators is independent of the choice of topology in $E'$, and in all the constructions described above the space $E'$ could be 
equipped not only with the strong topology $\beta(E', E)$, but also with the weak topology $\sigma(E', E)$ or any intermediate topology, for example, the Mackey topology 
$\tau(E', E)$. The class of Fredholm operators is independent of the choice of topology in $E$, and is denned only by the dual pair $(E', E)$; in this sense the class 
of Fredholm operators is ''more natural'' than the class of nuclear operators [{\bfseries 5}].

\section*{\small 3. Approximation-type properties of locally convex spaces}

{\bf 3.1.} We denote by $L(E)$ the algebra of all continuous linear operators in $E$ equipped with the topology of precompact convergence (i.e., uniform conver-gence on 
precompact sets). The space $E$ possesses the {\it approximation property} [{\bfseries 4}] if the continuous linear operators of finite rank\footnote{Below, unless stipulated otherwise, we will assume that all operators of finite rank under consideration are continuous linear operators, and the terms ''operator of finite rank'' and ''finite-dimensional operator'' will be used synonymously.}  form a dense subset in $L(E)$. If $E$ is a Banach space, it possesses the approximation property if and only if for an arbitrary Banach space $F$, every compact operator $F \to E$ may be approximated in norm by operators of finite rank. The approximation hypothesis disproved by Enflo [{\bfseries 2}] asserted that all locally convex (or at least all Banach) spaces possess the approximation property.

We will say that a space $E$ possesses the {\it weak approximation property} if the set of operators of finite rank is dense in $L(E)$ relative to the topology of 
uniform convergence on all absolutely convex compact sets in $E$. This definition is borrowed from [{\bfseries 6}]. It is clear that the weak approximation property 
follows from the approximation property, and that these properties are equivalent for complete (in particular, Banach) spaces.

A Banach space $E$ possesses the {\it bounded approximation property} if there exists a net of operators of finite rank $\{I_{\nu}\}$ that converges to the identity operator in the strong operator topology (i.e., in the topology of simple convergence), and furthermore, the norms $\|I_{\nu}\|$ are uniformly bounded (in $\nu$) [{\bfseries 2}]. A 
somewhat stronger metric approximation property for Banach spaces has been formulated and investigated by Grothendieck [{\bfseries 4}].

We may formulate the {\it bounded approximation property} for an arbitrary locally convex space $E$ in the following way. We will say that $E$ possesses this property 
if there exists a net of finite-dimensional operators that converges to the identity operator in the weak operator topology and is bounded in this topology. 
From Proposition 2 proved below, it follows that in Banach spaces our definition is equivalent to the \textquotedblleft classical\textquotedblright \  definition of the bounded approximation 
property presented above, and that in complete metrizable Hausdorff spaces it is equivalent to a definition formulated in [{\bfseries 1}] for the case of nuclear 
Frechet spaces.

The existence of a basis is a well-known approximation-type property closely connected to other approximation properties [{\bfseries 2}], [{\bfseries 4}]. A sequence $\{e_i\}$ of vectors in a locally convex space $E$ in which each element $x \in E$ has the unique decomposition

\begin{equation}
x = \sum_{i=1}^{\infty} \xi_i  e_i ,
\end{equation}
where the linear forms $x \mapsto \xi_i$ are continuous, is called a {\it Schauder basis} in the locally convex space $E$. If $E$ is a Banach space, the continuity of these forms follows automatically from the uniqueness of the decomposition (5). If the series (5) converges only in the weak topology $\sigma(E, E')$, the basis $\{e_i\}$ is 
said to be a {\it weak basis}.

We denote by $e'_i$ the element in $E'$ corresponding to the linear form $x \mapsto \xi_i $, so that $\xi_i = \langle e'_i, x \rangle $, and the decomposition (5) may be rewritten in the form

\begin{equation}
x = \sum_{i=1}^\infty \langle e'_i, x \rangle e_i .
\end{equation}
The sequence $\{e'_i\}$ forms a Schauder basis in the space $E'$ equipped with the weak topology $\sigma(E', E)$. This sequence is said to be {\it biorthogonal} to the basis 
$\{e_i\}$.
\\

{\bf 3.2.} Let us prove several assertions that establish the interrelation betwe-en the different approximation-type properties. Since every convergent coun-table 
sequence is bounded, the following assertion follows directly from the definition of the bounded approximation property.
\\

{\bf Proposition 1.} {\it If in a locally convex space $E$ there exists a countable sequence of finite-dimensional operators that converges to the identity ope-rator in 
the weak operator topology, $E$ possesses the bounded approximation property.}
\\

{\bf Corollary 1.} {\it If a locally convex space contains a Schauder basis (at least a weak Schauder basis), it possesses the bounded approximation property.}
\\

In fact, suppose that $\{e_i\}$ is a Schauder basis in $E$, and let $\{e'_i\}$ be the biorthogonal sequence in $E'$. We set $I_n \colon\  x \mapsto \sum_{i=1}^n  \langle e'_i, x\rangle e_i$. Equality (6) asserts that the sequence $\{I_n\}$ of finite-dimensional operators converges to the identity operator in the weak operator topology if the series (6) converges at least weakly. Now it is only necessary to refer to Proposition 1 to complete the proof.
\\

{\bf Proposition 2.} {\it $A$ space $E$ possesses the bounded approximation property if and only if there exists a net of operators of finite rank in $E$ that converges to 
the identity operator in the strong operator topology (i.e., in the topology of simple convergence) and is bounded in this topology.}
\\

{\bfseries Proof.} The space $L(E)$ has the same dual space $E' \otimes E$ relative to the strong and weak operator topologies (see [{\bfseries 12}], Ch. IV, sec. 4). Therefore, bounded sets and the closures of convex sets coincide for both topologies; this follows from the Hahn-Banach and Mackey-Arens theorems (see [{\bfseries 12}], Ch. IV, sec. 3). Suppose that $\{I_{\nu}\}$ is a net of finite-dimensional operators that converges weakly to the identity operator and is bounded. It is clear that the convex hull of this set is bounded, consists of operators of finite rank, and contains the identity operator in its closure (weak and strong). From the bounded approximation 
property, therefore, it follows that there exists a net of finite-dimensional operators that converges to the identity operator in the strong operator 
topology and is bounded in this topology. The converse assertion is self-evident.
\\

{\bf Proposition 3.} {\it If the space $E$ is barrelled and possesses the bounded appro-ximation property, it also possesses the approximation property.}
\\

{\bfseries Proof.} From Proposition 2 it follows that there exists a net $\{I_{\nu}\}$ of finite-dimensional operators in $E$ that converges to the identity operator in the 
topology of simple convergence and is bounded in this topology. Since $E$ is a barrelled space, it follows from the fact that the net is bounded that the set 
$\{I_{\nu}\}$ is equicontinuous. The topologies of simple and precompact convergence coincide on any equicontinuous subset in $L(E)$ (Banach-Steinhaus principle, see [{\bfseries 12}], Ch. lll, sec. 4). Therefore, the net $\{I_{\nu}\}$ converges to the identity operator in the topology of precompact convergence. Hence, it easily follows that for any operator $A \in L(E)$ the net $\{AI_{\nu}\}$ converges to the operator $A$ in $L(E)$; since the $AI_{\nu}$ are finite-dimensional operators, the approximation property is 
satisfied.
\\

{\bf Proposition 4.} {\it If the space $E$ is complete, metrizable, and separable, it possesses the bounded approximation property if and only if there exists a countable 
sequence of operators of finite rank that converges to the identity operator in the strong operator topology in $L(E)$ (and, automatically, in the topology of 
precompact convergence).}
\\

{\bfseries Proof.} If a sequence of the type specified in Proposition 4 exists, the boun-ded approximation property is satisfied by virtue of Proposition 1. Since any 
complete metrizable space is a barrelled space, it follows from the bounded approximation property, Proposition 2, and the Banach-Steinhaus principle that 
there exists an equicontinuous set $H$ in $L(E)$ that consists of operators of finite rank and that contains the identity operator in its closure relative to the 
topology of simple convergence. The restriction of this topology to $H$ is metrizable, i.e., coincides with the restriction to $H$ of the topology of metrizable 
simple convergence on a countable total subset in $E$. Moreover, the restrictions to $H$ of the topologies of simple and precompact convergence coincide (see 
[{\bfseries 12}], Ch. lll, props. 4.5 and 4.7). Since $H$ is metrizable, there exists a countable sequence of elements from $H$ that converges to the identity operator in 
the topology of simple and precompact convergence, as claimed.
\\

{\bf 3.3.}    The next lemma is due to Grothendieck.
\\

{\bf Lemma 1.} {\it Suppose that the space $E$ possesses a basis of absolutely convex neighborhoods of zero such that for any neighborhood U from this basis the Banach 
space $\widehat{E_U}$ possesses the approximation property. Then $E$ possesses the approximation property.}
\\

The proof of this lemma may be found in [{\bfseries 4}], Ch. 1, sec. 5.1, and in [{\bfseries 12}], Ch. lll, sec. 9.
\\

{\bf Proposition 5.} {\it If the topology of a locally convex space $E$ coincides with the weak topology $\sigma(E, E')$, this space possesses the approximation property.}
\\

{\bfseries Proof.} It is easily seen that for any neighborhood $U$ of zero in $E$, the space $\widehat{E_U}$ is finite-dimensional and, consequently, possesses the approximation pro-perty. Therefore, Proposition 5 follows from Lemma 1.

\section*{\small 4. Principal results}

{\bf 4.1.} The space $E' \bar{\otimes} E$ is a two-sided $S(E)$-module. The product $Au$, where $u \in E' \bar{\otimes} E$ and $A \in S(E)$, is defined so that the continuous linear mapping $u \mapsto Au$ corresponds to a separately continuous bilinear mapping $E' \times E \to E' \bar{\otimes} E$ under which the pair $(x', x)$ is carried into the element $x' \otimes Ax$ (cf. above, \S2.4 and \S2.5). Analogously, the product $uA$ is defined by means of the bilinear mapping $(x', x) \mapsto {^{t}A}x' \otimes x$, where $ ^{t}A$ is the operator adjoint (dual) to $A$. It is easily seen that, for any Fredholm kernel $u = \sum \lambda_i  x'_i  \otimes x_i$ , we have

$$
AÀu = \sum \lambda_i x'_i \otimes Ax_i \qquad \mbox {and} \qquad uA = \sum \lambda_i ~ {^{t}A}x'_i \otimes x_i
$$

Hence, it follows that if $u \in E' \overbracket[1pt][1pt]{\otimes} E$, then $Au \in E' \overbracket[1pt][1pt]{\otimes} E$ and $uA \in E' \overbracket[1pt][1pt]{\otimes} E$. From the definition of the trace tr$u$ of the tensor $u \in E' \bar{\otimes} E$ and the mapping $\Gamma \colon\  E' \overbracket[1pt][1pt]{\otimes} E \to S(E)$ (cf. above, \S2), the validity of the following relations may be derived directly:

\begin{equation}
\mbox{tr}(uA) = \mbox{tr}(Au) = \sum_{i=1}^{\infty} \lambda_i  \langle x'_i,  Ax_i \rangle,
\end{equation}

\begin{equation}
\Gamma(uA) = \Gamma(u)A; \qquad  \Gamma(Au) = A\Gamma(u),
\end{equation}
where $A \in S(E)$ and $u$ is a Fredholm kernel with a decomposition

$$
u =  \sum_{i=1}^{\infty} \lambda_i x'_i \otimes x_i
$$
of the form of (4).

The semi-norms $A \mapsto \vert$tr$(Au)\vert$ define the weak operator topology in $S(E)$ if the element $u$ runs through the algebraic tensor product $E' \otimes E$. Consider the topology in $S(E)$ specified by the semi-norms of this type, where $u$ runs through the space $E' \overbracket[1pt][1pt]{\otimes} E$ of Fredholm kernels. As in [{\bfseries 7}], this topology will be called the {\it ultra-weak topology}. In the case in which {\it E} is a Hilbert space, this topology coincides with the ultra-weak topology in the sense of the theory of von 
Neumann algebras.

The next assertion may be found in [{\bfseries 4}], [{\bfseries 7}].
\\

{\bf Proposition 6.} {\it Every Fredholm operator in the space $E$ possesses a well-defined trace if and only if the set of operators of finite rank is dense in $S(E)$ 
relative to the ultra-weak topology.}
\\

For the sake of completeness, we present the proof of this assertion. We will need the following lemmas.
\\

{\bf Lemma 2.} {\it Suppose that $F$ is an operator of finite rank in $E$, and that $F \colon \  x \mapsto \sum_{j= 1}^m \langle y'_j, x \rangle y_j$, where $y_j \in E, y'_j \in E'$. Then, for any Fredholm kernel $u \in E' \overbracket[1pt][1pt]{\otimes} E$, the relation} tr{\it $(Fu) = \sum_{j=1}^m  \langle y'_j, \Gamma(u)y_j \rangle = $} sp{\it $(F\Gamma(u))$ is satisfied.}
\\

{\bf Lemma 3.} {\it The following assertions are equivalent: (1) $\Gamma(u) = 0$; \\(2)} tr{\it $(uF) = 0$ for any operator $F$ of finite rank in $S(E)$.}
\\

{\bfseries Proof of Lemma 2.} Suppose that the Fredholm kernel $u$ can be represe-nted in the form (4), and that the operator $A = \Gamma(u)$ has the corresponding decomposition 
(1). Since the functionals $y'_j$ are continuous, it follows from (1) that

$$
\langle y'_j, Ay_j \rangle = \sum_{i=1}^{\infty} \lambda_i \langle x'_i, y_j\rangle \langle y'_j, x_i\rangle.
$$

Hence

$$
\sum_{j=1}^m \langle y'_j, \Gamma(u)y_j \rangle = \sum_{j=1}^m \sum_{i=1}^{\infty} \lambda_i \langle x'_i, y_j\rangle \langle y'_j, x_i\rangle.
$$

On the other hand,

$$
\mbox{tr}(Fu) = \sum_{i=1}^{\infty} \lambda_i \langle x'_i, Fx_i \rangle = \sum_{i=1}^{\infty} \sum_{j=1}^m  \lambda_i \langle x'_i, y_j\rangle \langle y'_j, x_i\rangle.
$$

Since the index $j$ runs through a finite set of values, summation with respect to $i$ may be interchanged with summation with respect to $j$, so that it only 
remains for us to note that the quantity $\sum_{j=1}^m \langle y'_j, \Gamma(u)y_j \rangle$ coincides with the ordinary trace sp$(F(u))$  of  the  finite-dimensional operator  $F\Gamma(u)$.
\\

{\bfseries Proof of Lemma 3.} If $\GammaÃ(u) = 0$, then by virtue of Lemma 2, tr$(uF) =$ tr$(Fu) =$ sp$(F \cdot \Gamma(u)) = 0$ for any finite-dimensional operator $F$. To prove the converse assertion, set $F \colon\  x \mapsto \langle e', x\rangle e$, where $e \in E$ and $e' \in E'$. In this case, tr$(uF) = \langle e', \Gamma(u)å)\rangle $. If tr$(uF) = 0$ for all operators $F$ of rank 1, then all the matrix elements $\langle e', \Gamma(u)e\rangle $ of the operator $\Gamma(u)$ are equal to 0, i.e., $\Gamma(u) = 0$, as claimed.
\\

{\bfseries Proof of Proposition 6.} From the definition of the ultra-weak topology, it follows that any linear functional $l$ on $S(E)$ that is continuous in this topology 
has the form $l(A) =$ tr$(uA)$, where $u \in E' \overbracket[1pt][1pt]{\otimes} E$. Therefore, if the operators of finite rank are not dense in $S(E)$ relative to the ultra-weak topology, by the Hahn--Banach theorem a continuous linear functional $l(A) =$ tr$(uA)$ may be found such that $l(A) = 0$ for all operators of finite rank and $l(1) \ne 0$, where 1 is the identity operator. It is clear that the Fredholm kernel $u$ has a nonzero trace, tr$u =$ tr$(u \cdot 1) = l(1) \ne 0$, but $\Gamma(u) = 0$ by virtue of Lemma 3. Thus, the 
trace cannot be well-defined. But if the operators of finite rank are dense in $S(E)$ relative to the ultra-weak topology, a continuous functional that 
vanishes on this dense set vanishes automatically on the identity operator as well, i.e., it follows from the condition $\Gamma(u) = 0$ that (by virtue of Lemma 3) 
tr $u = 0$. In this case, the trace is well-defined (cf. \S2.5), as claimed.
\\

{\bfseries Remark.} Proposition 6 may be generalized in the following way, as discussed in [{\bfseries 10}]. Let $T$ be a one-sided $S(E)$-submodule in $E' \bar{\otimes} E$. The topology in $S(E)$ and $L(E)$ specified by the semi-norms $A \mapsto \vert$tr$(uA)\vert $, where $u \in T$, is said to be  {\it  T-weak}.  Suppose that the canonical mapping $\Gamma \colon\  E' \otimes E \to S(E)$ (cf. \S2) may be continued to the subspace $T \supset E' \otimes E$. In this case, all the operators in $\Gamma(T)$ possess a well-defined trace if and only if the identity operator in $E$ may be approximated by operators of finite rank in the $T$-weak topology.

Thus, one can easily find approximation properties that provide an ade-quate description of the conditions under which it is possible for the trace to be 
well-defined for different classes of operators. However, the real problem is to find conditions that not only are of theoretical interest but also may be 
effectively verified.
\\
\\
\\

{\bf 4.2.}
\\

{\bf Lemma 4.} {\it Suppose that the Fredholm kernel è can be represented in the form (4) so that the sequence $\{x'_i\}$ is equicontinuous and $\Gamma(u)$ is a nuclear operator. Suppose that $\{A_{\nu}\}$ is a net in $S(E)$ that converges to the operator $A$ in the topology of uniform convergence on all absolutely convex compact sets. Then} tr{\it $(uA_{\nu) }\to$} tr{\it $(uA)$.}
\\

{\bfseries Proof.} We select a decomposition of the form (4) for the tensor $u$ such that the sequence $\{x_i\}$ is contained in some compact absolutely convex set $K$ in $E$, and such that the sequence $\{x'_i\}$ is contained in the polar $U^0$ of some absolutely convex neighborhood $U$ of zero in $E$. This decomposition always exists for the tensor $u$ (cf. [{\bfseries 4}], Ch. 1, sec. 3.2). The neighborhood $U$ absorbs the compact set $K$, i.e., there is a number $\mu > 0$ such that $Ax_i \in \mu U$. Since  the net $\{A_{\nu\}}$ converges uniformly on $K$, there exists an index ${\nu}_0$ such that $(A_{\nu} - A)x_i \in \mu U$ for $\nu \geq \nu_0$ for all values of the index $i$. Therefore, $A_{\nu} x_i \in 2\mu U$ and, consequently, $\vert \langle x'_i, A_{\nu} x_i\rangle \vert \leq 2\mu$ for all $i$ and all $\nu \geq \nu_0$, i.e., the quantity $\langle x'_i, A_{\nu} x_i \rangle $ is uniformly bounded, so that by Lebesgue's well-known theorem, it is possible to pass to the limit under the summation $\sum_{i=1}^{\infty} \lambda_i \langle x'_i, A_{\nu} x_i \rangle $, and this limit is equal to $\sum_{i=1}^{\infty} \lambda_i \langle x'_i, Ax_i \rangle $. In other words, tr$(uA_{\nu}) \to$ tr$(uA)$ by virtue of (7), as claimed.
\\

{\bf Theorem 1.} {\it Suppose that the space $E$ possesses the weak approximation property, i.e., there exists a net of operators of finite rank that converges to the 
identity operator in the topology of uniform convergence on all absolutely convex compact sets in $E$. Then every nuclear operator in $E$ has a well-defined 
trace. If $\{B_{\nu}\}$ is a net in $L(E)$ that converges to the operator $B \in L(E)$ in this topology, then} tr{\it $(AB_{\nu}) \to$} tr{\it $(AB)$ for any nuclear operator $A$.}
\\

{\bfseries Proof.} Suppose that $\{I_{\nu}\}$ is a net of finite-dimensional operators that converges uniformly to the identity operator on all absolutely convex compact sets, and $A = \Gamma(u)$ is an arbitrary nuclear operator, where the tensor $u$ satisfies the premise of Lemma 4. By virtue of this lemma and Lemma 2,

$$
{\mbox{tr}}A =\sum_{i=1}^{\infty} \lambda_i \langle x'_i , x_i \rangle =\lim_{\nu} \sum_{i=1}^{\infty} \lambda_i \langle x'_i , I_{\nu} x_i \rangle = \lim_{\nu} {\mbox tr}(uI_{\nu}) = \lim_{\nu} \mbox{sp}(AI_{\nu}).
$$

Since the quantities sp$(AI_{\nu})$ and lim sp$(AI_{\nu})$ are independent of the tensor $u$ and of decompositions of the form (1) and (4), the quantity tr$A$ is well-defined. The second assertion of the theorem also follows from Lemma 4.
\\

{\bf Corollary 2.} {\it If the space $E$ possesses the approximation property, every nuclear operator in $E$ possesses a well-defined trace.}
\\

{\bf 4.3.} Let us demonstrate that the approximation properties do not suffice for ensuring that every Fredholm operator possesses a well-defined trace. Consider 
the reflexive Banach space $E$ constructed by Enflo [{\bfseries 2}]; this space does not possess the approximation property. We denote by $\widetilde{E}$ the same space, but now equipped with the weak topology $\sigma (E, E')$. The space $\widetilde{E}$ is quasi-complete (see [{\bfseries 12}], Ch. IV, prop. 5.5) and possesses the approximation property by virtue of Proposition 5. Note that the Fredholm operators in $E$ and $\widetilde{E}$ are the same (although, by the same token, all nuclear operators in $\widetilde{E}$ are finite-dimensional). The next assertion follows, therefore, from Grothen-dieck's results on the relation between the uniqueness problem and the appro-ximation problem for Banach 
spaces.
\\

{\bf Theorem 2.} {\it There exists a quasi-complete locally convex space with the approximation property in which the trace cannot be well-defined for all Fred-holm 
operators.}
\\

{\bf Theorem 3.} {\it Suppose that a locally convex space E possesses the bounded approximation property. Then the trace of every Fredholm operator in $E$ is well-defined. If the net $\{A_{\nu}\}$ in $S(E)$ converges to the operator $A$ in the weak operator topology and is bounded in this topology, $\lim_{\nu}$} tr{\it $(FA_{\nu}) =$} tr{\it $(FA)$ is satisfied by every Fredholm operator $F$.}
\\

Theorem 3 can be derived from Lemma 2 and the following lemma in precisely the same way as Theorem 1 was derived from Lemmas 2 and 4.
\\

{\bf Lemma 5.} {\it If the net $\{A_{\nu}\}$ in $S(E)$ converges to the operator $A$ in the weak operator topology and is bounded in this topology,} tr{\it $(uA_{\nu}) \to$} tr{\it $(uA)$ for any tensor  $u \in E' \overbracket[1pt][1pt]{\otimes} E$.}
\\

{\bf Corollary 3.} {\it If the net $\{A_{\nu}\}$ in $S(E)$ converges in the weak operator topology and is bounded in this topology, it converges in the ultra-weak topo-logy as well.}
\\

{\bf Corollary 4 [7].} {\it If the identity operator in $E$ is the limit in the weak operator topology of a countable sequence $\{I_n\}$ of operators of finite rank, every 
Fredholm operator $F$ in $E$ possesses a well-defined trace and}
$$
\mbox{tr} F = \lim_{n \to \infty} \mbox{tr} (F I_n) = \lim_{n \to \infty} \mbox{sp} (FI_n) .
$$

Lemma 5 is proved in precisely the same way as Lemma 4. The uniform boundedness of the quantities $\vert \langle x'_i, A_{\nu} x_i \rangle \vert $ is derived directly from the boun-dedness of the sets $\{A_{\nu}\}$, $\{x_i\}$ and $\{x'_i\}$. In fact, Lemma 5 and Corollary 3 were proved in [{\bfseries 7}] (cf. the proof of Proposition 4.2 and the remark on page 343 in [{\bfseries 7}]). Corollary 3 follows from Lemma 5 and the definition of the ultra-weak topology. Corollary 4 follows from Theorem 3 and Proposition 1 or from Proposition 6 and Corollary 3, since every countable sequence in $S(E)$ that is convergent in the weak topology is bounded, so that if the premise of Corollary 4 is satisfied, then so are the premises of Theorem 3. Note that Theorem 3 also follows from Proposition 6 and Corollary 3.
\\

{\bf Corollary 5.} {\it There exists a quasi-complete locally convex space that posses-ses the approximation property but does not possess the bounded appro-ximation 
property.}
\\

This corollary follows from Theorems 2 and 3. There even exists a Banach space with the approximation property that does not possess the bounded 
approximation property. This assertion has been proved in [{\bfseries 3}].
\\

{\bf 4.4.} Suppose that the space $E$ possesses a Schauder basis $\{e_i\}$ (at least a weak Schauder basis) and let $\{e'_i\}$ be the sequence in $E'$ that is biorthogonal to 
this basis. Set $I_n \colon\  x \mapsto \sum_{i=1}^n \langle e'_i, x \rangle e_i$. Clearly, $\{I_n\}$ is a sequence of operators of finite rank that converges to the identity operator in the weak operator topology. From Corollary 4, it follows that the trace of every Fredholm operator $F$ in $E$ is well-defined; moreover,

$$
\mbox{tr} F = \lim_{n \to \infty} \mbox{tr}(FI_n)= \lim_{n \to \infty}  \sum_{i=1}^n \langle e'_i, Fe_i \rangle = \sum_{i=1}^n \langle e'_i, Fe_i \rangle.
$$

The quantity

$$
\sum_{i=1}^{\infty} \langle e'_i, Fe_i \rangle
$$
is called the {\it matrix trace} of the operator $F$. Thus, we arrive at the following assertion.
\\

{\bf Corollary 6 [8].} {\it If the space $E$ possesses a Schauder basis (or at least a weak Schauder basis), then every Fredholm operator in $E$ possesses a well-defined trace 
that coincides with the matrix trace.}
\\

This assertion has been proved for absolute bases in Banach spaces by A. S. Markus and V. I. Matsaev [{\bfseries 11}] using a more complicated technique.

\renewcommand\refname{References}

\end{document}